
\input amssym.def
\input amssym
\input psfig
\magnification=1100
\baselineskip = 0.25truein
\lineskip = 0.01truein
\vsize = 8.5truein
\voffset = 0.2truein
\parskip = 0.10truein
\parindent = 0.3truein
\settabs 12 \columns
\hsize = 5.4truein
\hoffset = 0.4truein

\setbox\strutbox=\hbox{%
\vrule height .708\baselineskip
depth .292\baselineskip
width 0pt}
\font\caps=cmcsc10

\font\bigtenrm=cmr10 at 14pt
\font\bigmaths=cmmi10 at 16pt

\def\sqr#1#2{{\vcenter{\vbox{\hrule height.#2pt
\hbox{\vrule width.#2pt height#1pt \kern#1pt
\vrule width.#2pt}
\hrule height.#2pt}}}}
\def\square{\mathchoice\sqr46\sqr46\sqr{3.1}6\sqr{2.3}4}

\centerline{\bigtenrm LARGE GROUPS, \textfont1=\bigmaths
PROPERTY ($\tau$) }
\centerline{\bigtenrm AND THE HOMOLOGY GROWTH OF SUBGROUPS}
\tenrm
\vskip 14pt
\centerline{MARC LACKENBY}
\vskip 18pt

\centerline{\caps 1. Introduction}
\vskip 6pt

In this paper, we investigate the homology of finite index subgroups
$G_i$ of a given finitely presented group $G$. We fix a prime $p$,
denote the field of order $p$ by ${\Bbb F}_p$, and define
$d_p(G_i)$ to be the dimension of $H_1(G_i; {\Bbb F}_p)$.
We will be interested in the situation where $d_p(G_i)$ grows
fast as a function of the index $[G:G_i]$. Specifically,
we say that a collection of finite index subgroups 
$\{ G_i \}$ has {\sl linear growth of mod $p$ homology} if
$\inf_i d_p(G_i)/ [G:G_i]$ is positive.
A major class of groups $G$ having such a collection of
subgroups are those that are {\sl large}.
By definition, this means that $G$ has a
finite index subgroup that admits a surjective
homomorphism onto a free non-abelian group.
One might wonder whether
largeness is {\sl equivalent} to the existence
of some nested sequence of finite index subgroups $\{ G_i \}$
with linear growth of mod $p$ homology for some prime $p$.
We will show that this is true if
one is willing to make extra hypotheses.
Firstly, we suppose that each $G_{i+1}$ is normal in
$G_i$ and has index a power of $p$. Secondly, we use the notion
of Property $(\tau)$. This is an important
group-theoretic concept, first defined by Lubotzky and
Zimmer [9], with connections
to graph theory, representation theory and differential
geometry. We will recall its definition in Section 2.
We will show that the largeness
of a finitely presented group can be
characterised in terms of linear growth of mod $p$ homology
and the failure of Property $(\tau)$.
Our main theorem is the following.

\noindent {\bf Theorem 1.1.} {\sl Let $G$ be a finitely presented group,
let $p$ be a prime and suppose that $G \geq G_1 \triangleright G_2 \triangleright \dots$
is a nested sequence of finite index subgroups, such
that each $G_{i+1}$ is normal in $G_{i}$ and has index a power of $p$. 
Suppose that $\{ G_i \}$ has linear growth of mod $p$ homology. Then,
at least one of the following must hold:
\item{(i)} some $G_i$ admits a surjective homomorphism onto $({\Bbb Z}/p{\Bbb Z})
\ast ({\Bbb Z}/p{\Bbb Z})$ and some normal subgroup of $G_i$, with index
a power of $p$, admits a surjective homomorphism onto a non-abelian free group;
in particular, $G$ is large;
\item{(ii)} $G$ has Property $(\tau)$ with respect to $\{G_i \}$.

}

\vfill\eject
The two possible conclusions in this theorem can be viewed
as a `win/win' scenario. On the one hand, largeness is a very useful
property. For example, it implies that the group has super-exponential
subgroup growth and infinite virtual first Betti number. On the
other hand, Property $(\tau)$ has many interesting applications,
for example to spectral geometry and random walks [6].

As an almost immediate consequence of Theorem 1.1, we obtain
the following characterisation of large finitely
presented groups. We will give a proof of this, assuming Theorem 1.1,
in Section 2.

\noindent {\bf Theorem 1.2.} {\sl Let $G$ be a finitely presented group.
Then the following are equivalent:
\item{(i)} $G$ is large;
\item{(ii)} there exists a sequence of finite index subgroups,
$G \geq G_1 \triangleright G_2 \triangleright \dots,$ and a prime $p$ such that
\itemitem{(1)} $G_{i+1}$ is normal in $G_i$ and has index a power of 
$p$, for each $i$;
\itemitem{(2)} $G$ does not have Property $(\tau)$ with respect
to $\{G_i \}$; and
\itemitem{(3)} $\{ G_i \}$ has linear growth of mod $p$ homology.

}

Theorem 1.1 can also be used to provide a substantial
class of groups that have Property $(\tau)$ with respect
to some nested sequence of finite index subgroups.

\noindent {\bf Theorem 1.3.} {\sl Let $G$ be a finitely
presented group and let $p$ be a prime. Suppose
that $G$ has an infinite nested sequence of subnormal subgroups,
each with index a power of $p$, and with linear
growth of mod $p$ homology. Then $G$ has such a
sequence that also has Property $(\tau)$.}

Theorem 1.1 bears a strong resemblance to another result
of the author. In [4], the following was proved:

\noindent {\bf Theorem 1.4.} {\sl Let $G$ be a finitely
presented group, and let $\{ G_i \}$ be a nested sequence
of finite index normal subgroups. Then at least one of the
following holds:
\item{(i)} $G_i$ is an amalgamated free product
or HNN extension for all sufficiently large $i$;
\item{(ii)} $G$ has Property $(\tau)$ with respect
to $\{ G_i \}$;
\item{(iii)} $\inf_i d(G_i)/[G:G_i]$ is zero.}

Here, $d( \ )$ is the rank of a group, which is the
minimal size of a generating set. In this paper,
$d_p( \ )$ plays this r\^ole; using $d_p( \ )$ rather
than $d( \ )$, we strengthen
(i) to deduce that $G$ is large.
Not only are the statements of Theorems 1.1 and 1.4
very similar, but also their proofs follow similar lines, although
the proof of Theorem 1.1 is more complicated.
The geometry and topology of Schreier coset graphs play a central
r\^ole in both arguments.
The main difference is that a key application of
the Seifert - van Kampen theorem in the proof of the
Theorem 1.4 is replaced by the Mayer - Vietoris theorem
with mod $p$ coefficients in the proof of Theorem 1.1.

There is an interesting application of Theorem 1.1
to low-dimensional topology and geometry. A major
area of research in this field is the study of
lattices in ${\rm PSL}(2,{\Bbb C})$ (or, equivalently,
finite-volume hyperbolic 3-orbifolds). An important unsolved
problem asks whether any such lattice is
a large group. In [5], it was shown that
if such a lattice contains a torsion element
then it has a nested sequence $\{ G_i \}$ of finite index 
subgroups with linear growth of mod $p$ homology, for
some prime $p$. Moreover, 
these subgroups are all normal in $G_1$ and have
index a power of $p$. Thus, we deduce
from Theorem 1.1 that {\sl either} $G$ has
Property $(\tau)$ with respect to $\{ G_i \}$
{\sl or} that $G$ is large. In [5], we show that
the following conjecture of Lubotzky and Zelmanov,
which we have termed the GS-$\tau$ Conjecture,
implies that we can arrange that the former possibility
does not arise.

\noindent {\bf Conjecture 1.5.} (GS-$\tau$ Conjecture) {\sl Let $G$ be  
a group with finite presentation $\langle X | R \rangle$,
and let $p$ be a prime.
Suppose that $d_p(G)^2/4 > |R| - |X| + d_p(G)$.
Then $G$ does not have Property $(\tau)$
with respect to some infinite nested sequence $\{ G_i \}$
of normal subgroups with index a power of $p$.}

Thus, Theorem 1.1 and the argument in [5] give the
following result.

\noindent {\bf Theorem 1.6.} {\sl The GS-$\tau$
Conjecture implies that any lattice in 
${\rm PSL}(2, {\Bbb C})$ with torsion
is large.}

It is natural to ask which finitely generated groups
$G$ have a sequence of subnormal subgroups, each with
index a power of $p$ and with linear growth of mod $p$ homology.
We prove a stronger version of the following result in Section 8, 
which gives an alternative characterisation of these
groups.

\noindent {\bf Theorem 1.7.} {\sl Let $G$ be a finitely generated
group, and let $p$ be a prime. Then the following
are equivalent:
\item{(i)} $G$ has an infinite nested sequence of subnormal subgroups,
each with index a power of $p$, and with linear growth
of mod $p$ homology;
\item{(ii)} the pro-$p$ completion of $G$ has 
exponential subgroup growth.

}

Combining Theorems 1.3 and 1.7, we have the following interesting corollary.

\noindent {\bf Corollary 1.8.} {\sl Let $G$ be a finitely presented
group, and let $p$ be a prime. Suppose that the pro-$p$ completion of
$G$ has exponential subgroup growth. Then $G$ has a nested sequence
of subnormal subgroups, each with index a power of $p$, which
has Property $(\tau)$.}

Property $(\tau)$ plays a prominent r\^ole in the statement of
Theorem 1.1. But one might wonder to what extent it is needed.
Might it be true that conclusion (i) of Theorem 1.1 always holds?
We will see how this question relates to error-correcting codes.
We will show that if (i) does not hold, then an infinite collection
of linear codes can be constructed that are `asymptotically good'.
These are very important in the theory of error-correcting codes,
because they have large rate and large Hamming distance.
More details of this relationship can be found in Section 6.

We now briefly describe the plan of the paper.
In Section 2, we recall the definition of
Property $(\tau)$, and then go on to prove
Theorems 1.2 and 1.3 from Theorem 1.1.
In Section 3, we give a necessary and sufficient
topological condition on a finite connected 2-complex
(satisfying some generic conditions) for its fundamental
group to admit a surjective homomorphism onto a 
non-abelian free group. This is a key step in the
proof of Theorem 1.1, which is presented in Sections
4 and 5. Section 5 in particular is the heart of 
the paper. In Section 6, we establish a link
between large groups and error-correcting codes.
In Section 7, we show that the assumption of finite
presentability in Theorems 1.1 and 1.3 cannot be weakened to being finitely
generated. This is because the (generalised) lamplighter
group $({\Bbb Z}/p{\Bbb Z}) \wr {\Bbb Z}$, which is
finitely generated, satisfies the remaining hypotheses of
Theorem 1.1 and 1.3 but satisfies none of their conclusions.
Finally, in Section 8, we relate linear growth of
mod $p$ homology to the subgroup growth of the group's
pro-$p$ completion. 

I am grateful to Jim Howie and Alex Lubotzky who suggested to me the
examples in Section 7. I would also like to thank Andrei
Jaikin who suggested an improvement to an earlier version
of Proposition 4.2.

\vskip 18pt
\centerline{\caps 2. Property $(\tau)$}
\vskip 6pt

In this section, we recall the definition of
Property $(\tau)$, and then go on to deduce
Theorems 1.2 and 1.3 from Theorem 1.1.

Let $G$ be a finitely generated group, and let
$\{ G_i \}$ be a collection of finite index
subgroups. Let $S$ be a finite generating set
for $G$, and let $X(G/G_i; S)$ be the
Schreier coset graph for $G/G_i$ with
respect to $S$.

Property $(\tau)$ is defined in terms of
the geometry of these graphs. Specifically,
we will look at subsets $A$ of their vertex
set and consider $\partial A$, which is defined
to be the set of edges with one endpoint in $A$
and one not in $A$. The {\sl Cheeger constant} 
$h(X)$ of a finite graph $X$ is defined to be
$$h(X) = \min \left \{ {|\partial A| \over |A|}:
A \subset V(X) \hbox{ and } 0 < |A| \leq |V(X)|/2 \right \},$$
where $V(X)$ is the vertex set of $X$.
Then $G$ is said to have {\sl Property $(\tau)$ with
respect to $\{ G_i \}$} if $\inf_i h(X(G/G_i; S))$ is
strictly positive,
for some finite generating set $S$ for $G$. It turns
out that if this holds for some finite generating
set then it holds for any finite generating set
(see Lemma 2.3 in [3] for example).

A basic example is the group $G = {\Bbb Z}$ and
its subgroups $G_n = n {\Bbb Z}$. Let $S = \{ 1 \}$.
Then $X(G/G_n; S)$ is a circular graph with
$n$ vertices and $n$ edges. It is clear that
$h(X(G/G_n;S)) \rightarrow 0$. Hence, $G$ does not
have Property $(\tau)$ with respect to $\{ G_n \}$.
In fact, $G$ does not have Property $(\tau)$ with
respect to any infinite subcollection of $\{ G_n \}$.

The following two lemmas are elementary and well known.

\noindent {\bf Lemma 2.1.} {\sl Let $G$ and $K$ be finitely
generated groups, and let $\phi \colon G \rightarrow K$
be a surjective homomorphism. Let $\{K_i \}$ be
a collection of finite index subgroups
of $K$. Then $K$ has Property $(\tau)$ with
respect to $\{ K_i \}$ if and only if
$G$ has Property $(\tau)$ with respect to
$\{ \phi^{-1}(K_i) \}$.}

\noindent {\sl Proof.} Let $S$ be a finite generating
set for $G$. Then $\phi(S)$ forms a finite generating
set for $K$. Now, $\phi$ induces a bijection
between the right cosets $G/\phi^{-1}(K_i)$ and $K/K_i$. This respects
right multiplication by elements of $G$.
Hence, the coset graphs $X(G/\phi^{-1}(K_i); S)$ and
$X(K/K_i; \phi(S))$ are isomorphic. The lemma follows
immediately. $\square$

\noindent {\bf Lemma 2.2.} {\sl Let $G$ be a 
finitely generated group, and let $K$ be
a finite index subgroup. Let $\{ K_i \}$ be
a collection of finite index subgroups of $K$.
Then $G$ has Property $(\tau)$ with respect to
$\{ K_i \}$ if and only if $K$ has Property
$(\tau)$ with respect to $\{ K_i \}$.}

\noindent {\sl Proof.} This is essentially contained in the proof
of Lemma 2.5 in [3], but we include the proof here for the
sake of completeness, and because we are explicitly dealing
here with subgroups that need not be normal.

Let $S$ be a finite generating set for $G$. Let $T$ be a 
maximal tree in $X(G/K; S)$. Then the edges not in $T$
form a finite generating set $\tilde S$ for $K$,
by the Reidermeister-Schreier process. For any subgroup
$K_i$ of $K$, $X(G/K_i; S)$ is a covering space of
$X(G/K;S)$. The inverse image of $T$ in $X(G/K_i; S)$
is a forest $F$. If one were to collapse each component of
this forest to a point, one would obtain $X(K/K_i; \tilde S)$.

Let $A$ be any non-empty subset of the vertex set of $X(K/K_i; \tilde S)$.
Its inverse image $\tilde A$ in $X(G/K_i;S)$ is a union
of components of $F$. It is clear that $|\tilde A| = {[G:K]|A|}$
and $|\partial \tilde A| = |\partial A|$. Hence,
$h(X(G/K_i; S)) \leq h(X(K/K_i; \tilde S))/[G:K]$.
So if $h(X(K/K_i; \tilde S))$ has zero infimum, then so does
$h(X(G/K_i;S))$.

Now consider a non-empty subset $B$ of the vertex set of $X(G/K_i; S)$
such that $|\partial B|/|B| = h(X(G/K_i; S))$ and 
$|B| \leq |V(X(G/K_i; S))|/2$.
Let $\overline B$ be the vertices of 
$X(G/K_i; S)$ lying in the union of those components of $F$
that intersect $B$. Thus, $\overline B$ clearly contains
$B$. If a component of $F$ lies in $\overline B$ but does not
lie entirely in $B$, then it contains an edge of $\partial B$.
Hence,
$$|B| \leq |\overline B| \leq |B| + [G:K]|\partial B|.$$
If an edge lies in $\partial \overline B$ but not in
$\partial B$, then it joins two different components
of $F$, at least one of which contains an edge of
$\partial B$. There are at most $2 |S| [G:K]$ edges
with an endpoint in this component of $F$. Hence,
$$|\partial \overline B| \leq (2 |S| [G:K] + 1) |\partial B|.$$
Now, $\overline B$ projects to a set of vertices in
$X(K/K_i; \tilde S)$ with size that has been scaled by a 
factor of $[G:K]^{-1}$ and with the same size boundary.
Hence,

$$\eqalign{
& h(X(K/K_i; \tilde S)) \cr
&\qquad \leq {|\partial \overline B|
\over [G:K]^{-1} \min \{ |\overline B|, |\overline B^c| \}} \cr
&\qquad \leq [G:K] (2 |S| [G:K] + 1) { |\partial B| \over
\min \{ |B| , |B^c| - [G:K]|\partial B| \} }\cr
&\qquad \leq [G:K] (2 |S| [G:K] + 1) \max \left \{ h, { h \over 1 - [G:K] h} \right \},\cr}$$
where $h = h(X(G/K_i; S))$, and provided that
$|B^c| - [G:K]|\partial B| > 0$. This assumption certainly
holds if $h < [G:K]^{-1}$. So, if $h(X(G/K_i; S))$ has zero
infimum, then so does $h(X(K/K_i; \tilde S))$. $\square$

We are now in a position to prove Theorems 1.2 and 1.3.

\noindent {\sl Proof of Theorem 1.2.} (ii) $\Rightarrow$ (i) is an immediate consequence
of Theorem 1.1. In the other direction, suppose that some finite index
subgroup $G_1$ of $G$ admits a surjective homomorphism $\phi_1$
onto a non-abelian free group $F$.
Let $\phi_2 \colon F \rightarrow {\Bbb Z}$ 
be projection onto the first free summand. 
Now, ${\Bbb Z}$ does not have Property $(\tau)$
with respect to $\{ p^i {\Bbb Z} \}$, by the
earlier example. Let $G_i$ be $\phi_1^{-1} \phi_2^{-1}(p^{i-1} {\Bbb Z})$.
Then, for each $i$, $G_{i+1}$ is normal in $G_i$ and has
index $p$. By Lemma 2.1,
$G_1$ does not have Property $(\tau)$ with
respect to $\{ G_i \}$. By Lemma 2.2,
$G$ also does not have Property $(\tau)$ with
respect to $\{ G_i \}$.
Now, $\phi_2^{-1}(p^{i-1} {\Bbb Z})$ forms a nested
sequence of finite index subgroups in $F$, and any such sequence
has linear growth of mod $p$ homology.
As each $G_i$ surjects onto $\phi_2^{-1}(p^{i-1} {\Bbb Z})$,
$d_p(G_i) \geq d_p(\phi_2^{-1}(p^{i-1} {\Bbb Z}))$. Hence, 
$\{ G_i \}$ has linear growth of mod $p$ homology.
$\square$

\noindent {\sl Proof of Theorem 1.3.} If the given sequence of
subgroups has Property $(\tau)$, we are
done. If not, then Theorem 1.1 implies that some finite
index subnormal subgroup of $G$, with index a power of $p$, admits a surjective
homomorphism onto a non-abelian free group $F$.
By passing to a smaller subgroup of $G$ if necessary,
we may assume that $F$ has arbitrarily large rank. 
We claim that $F$ then has a sequence of normal subgroups,
each with index a power of $p$, with linear growth of
mod $p$ homology and with Property $(\tau)$. Their inverse
images in $G$ form the required subgroups by Lemmas 2.1 and 2.2.
There are many ways to prove this claim. One is to use
the fact that ${\rm SL}(3, {\Bbb Z})$ has Property $(\tau)$
with respect to its principal congruence subgroups [6].
Let $K_n$ denote the level $p^n$ principal congruence
subgroup. Then $K_{n+1}$ is normal in $K_n$
and has index a power of $p$, for all $n \geq 1$.
If the rank of $F$ is large enough, it admits a surjective
homomorphism onto $K_1$. The inverse images of $K_n$ in $F$
then form the required subgroups. $\square$

\vskip 18pt
\centerline{\caps 3. Cocycles and large groups}
\vskip 6pt

In this section, we will study connected finite 2-complexes $K$
and give a necessary and sufficient topological condition
for $\pi_1(K)$ to admit a free non-abelian quotient.
{\sl We make convention throughout this paper} that the attaching map
of each 2-cell of $K$ is cellular; that is, the boundary
path of the 2-cell can be expressed as a concatenation of
a finite sequence of paths, each of which is a homeomorphism
onto a 1-cell of $K$. 

The necessary and sufficient
condition will be phrased in terms of {\sl regular
cocycles}. These are particularly nice representatives
of elements of $H^1(K)$. We will show that any such cohomology
class is represented by a regular cocycle.

A regular cocycle is just a non-empty finite graph $\Gamma$ embedded within $K$ in
a certain way, together with orientation information.
The graph must satisfy the following conditions:
\item{(i)} $\Gamma$ is disjoint from the 0-skeleton of $K$;
\item{(ii)} its vertices $V(\Gamma)$ are the intersection of
$\Gamma$ with the 1-skeleton of $K$;
\item{(iii)} for any 2-cell with quotient map
$i \colon D \rightarrow K$, where $D$ is a 2-disc, $D \cap i^{-1}(\Gamma)$
is a finite collection of properly embedded disjoint arcs with endpoints
precisely $\partial D \cap i^{-1}(V(\Gamma))$.

\noindent We then say that the graph is {\sl regularly embedded}.
A {\sl regular cocycle} is a regularly embedded graph with a transverse
orientation assigned to each arc in each 2-cell, with
the requirement that near each vertex of $\Gamma$, these
transverse orientations all coincide.

A regular cocycle determines an element of $H^1(K)$, as follows.
It assigns to each oriented 1-cell of $K$ a weight, which is
just its signed intersection number with $\Gamma$. The
total weight of the boundary of any 2-cell is clearly zero.
This therefore gives a well-defined cellular cocycle and
hence an element of $H^1(K)$.

Conversely, one may construct
a representative regular cocycle
for any element of $H^1(K)$, as follows. Pick a cellular cocycle
representing the cohomology class. This is just an assignment
of an integer weight to each oriented 1-cell, with the
property that the weights of the boundary of any 2-cell
sum to zero. For any 1-cell $e$, with weight $w(e)$, say,
place $|w(e)|$ vertices of $\Gamma$ on the interior of
$e$. Give $e$ an orientation, so that its weight is
non-negative. Assign the same transverse orientation 
to the vertices on $e$. Since the total evaluation
around each 2-cell is zero, there is a way to insert the
transversely oriented edges of $\Gamma$ into the 2-cells,
forming a regular cocycle.

Note that a connected regular cocycle represents
a non-trivial element of $H^1(K)$ if and only if it
is non-separating. For, if it is separating, then its
evaluation of any closed loop in $K$ is zero, and hence
it represents the trivial cohomology class. Conversely,
if it is non-separating, then its evaluation on some
closed loop is non-zero, and so the associated cohomology
class is non-trivial.

We say that a point $x$ in $K$ is {\sl locally
separating} if it has a connected neighbourhood
$U$ such that $U - x$ is disconnected.
The {\sl valence} of a 1-cell of $K$ is the
total number of times the 2-cells of $K$ run over
it. In the second half of the following result, we consider only
finite 2-complexes with no locally separating points
and no 1-cells with valence 1. Note that any finite
2-complex can be transformed into a finite 2-complex
with these properties, without changing
its fundamental group. For, we may replace each 0-cell with a 2-sphere and each
1-cell with a tube. Thus, any finitely presented
group arises as the fundamental group of
a finite 2-complex with these properties.

For a group $G$ and positive integer $n$, let
$\ast^n G$ denote the free product of $n$
copies of $G$. For a space $X$ with a basepoint,
let $\bigvee^n X$ denote the wedge of $n$ copies of 
$X$ glued along their basepoints.

\vfill\eject
\noindent {\bf Theorem 3.1.} {\sl 
Let $K$ be a finite connected 2-complex. Then $\pi_1(K)$
admits a surjective homomorphism onto $\ast^n {\Bbb Z}$
if $K$ contains $n$ disjoint regular cocycles whose union
is non-separating. Furthermore, the converse also holds, provided
$K$ has no
locally separating points and no 1-cells with valence 1.}

\noindent {\sl Proof.} Suppose first that $K$ contains $n$ disjoint
regular cocycles $C_1, \dots, C_n$ whose union is non-separating.
These have disjoint product
neighbourhoods $C_i \times [-1,1]$. Define a map $f \colon
K \rightarrow \bigvee^n S^1$, as follows. Away from 
$\bigcup (C_i \times [-1,1])$, send everything to the
central vertex of $\bigvee^n S^1$. On $C_i \times [-1,1]$,
first project onto the second factor $[-1,1]$, and then compose
this with the quotient map $[-1,1] \rightarrow S^1$ that identifies
the endpoints of the interval, and then
map this to the $i^{\rm th}$ circle of $\bigvee^n S^1$.
Pick a basepoint $b$ for $K$ away from the neighbourhoods
of the cocycles. We claim that the induced map 
$f_\ast \colon \pi_1(K,b) \rightarrow \ast^n {\Bbb Z}$
is a surjection. This is because the $i^{\rm th}$ free generator
of $\ast^n {\Bbb Z}$ may be realised by a loop that
starts at $b$, runs to $C_i$, crosses it transversely, and returns
to $b$. We may ensure that this is the only point of
intersection between the loop and $\bigcup C_i$, by
the hypothesis that $\bigcup C_i$ is non-separating.

Conversely, suppose that $\pi_1(K)$ admits a surjective
homomorphism onto $\ast^n {\Bbb Z}$. We will show that
this is induced by a map $f \colon K \rightarrow \bigvee^n S^1$.
Pick a basepoint $b$ for $K$ in the 0-skeleton.
Pick a maximal tree $T$ in the 1-skeleton of $K$. Let $f$ send
this tree to the central vertex of $\bigvee^n S^1$.  Each remaining
of edge $e$ of $K$, when oriented, determines an element of $\pi_1(K, b)$,
given by the path that starts at $b$, runs along $T$ to the initial vertex of
$e$, then along $e$, then back to $b$ by a path in $T$.
The image of this element of $\pi_1(K, b)$ under the given
homomorphism is an element of $\ast^n {\Bbb Z}$, which we may take to be a
reduced word. This then gives a path in $\bigvee^n S^1$. Define
the restriction of $f$ to $e$ to be this path. Since we started with a 
homomorphism $\pi_1(K) \rightarrow \ast^n {\Bbb Z}$,
the boundary of each 2-cell is sent a homotopically
trivial loop in $\bigvee^n S^1$, and hence, there is
a way to extend $f$ over the 2-cells. Pick points $p_1, \dots, p_n$,
one in each circle of $\bigvee^n S^1$, disjoint
from the central vertex. Then it is clear that we may
ensure that, for each $i$, $f^{-1}(p_i)$ is a regularly embedded
graph. Moreover, if we impose orientations on the circles,
then these graphs inherit transverse orientations, making
them regular cocycles $C_1, \dots, C_n$, say.
These cocycles are clearly disjoint, but their union
may not yet be non-separating.
The aim now is to modify $f$ by a homotopy,
thereby changing the cocycles $C_i$, to
ensure that this is the case.

Define a graph $Y$, whose vertices correspond to the components of the complement of $\bigcup C_i$.
Let its edges be in one-one correspondence
with the components of $\bigcup C_i$, and
where incidence between edges and vertices
in $Y$ is defined by topological incidence in $K$.
The edges inherit an orientation from $\bigcup C_i$, 
and also inherit a label $i$. We will modify $f$, thereby
giving new regular cocycles $C_i$, and
hence a new graph $Y$. At each stage, the number of
components of $\bigcup C_i$ will decrease, and so
this process is guaranteed to terminate. The aim
is to ensure that $Y$ satisfies the following condition:
\item{($\ast$)} no vertex of $Y$ has two edges pointing
into it with the same label, or two edges
pointing out of it with the same label.

Suppose now that ($\ast$) is violated. Let $E_1$ and $E_2$
be distinct components of $C_i$, say, both pointing
into the same component $X$ of $K - \bigcup C_i$.
Since $K$ contains no locally separating points, each
1-cell of $K$ has non-zero valence. Hence, neither $E_1$
nor $E_2$ is a point.
Pick an embedded arc $\alpha$, with one endpoint on
$E_1$ and the other endpoint on $E_2$, and with
interior in $X$. Since every 1-cell of
$K$ has valence at least two, every vertex of the
graphs $E_1$ and $E_2$ has valence at least two.
So neither graph is a tree.
Hence, each contains a point in the interior of
an edge such that removing that point from $E_i$
does not disconnect $E_i$. We may assume that
the endpoints of $\alpha$ are these two points.
Because $K$ has no locally separating
points, we may arrange for $\alpha$ to miss the 0-cells of $K$.
We may ensure that $\alpha$ intersects each 1-cell
in a finite collection of points, and each 2-cell
in a finite collection of arcs, each of which is properly
embedded, except the arc(s) containing the endpoints
of $\alpha$. Let $\alpha \times [-1,1]$ be a thickening
of $\alpha$, so that $(\alpha \times [-1,1]) \cap \bigcup C_i
= \partial \alpha \times [-1,1]$. We now modify $f$,
leaving it unchanged away from a small regular
neighbourhood of $\alpha \times [-1,1]$. In $\alpha \times [-1,1]$,
modify $f$ so that the intersection of the new $C_i$ with
$\alpha \times [-1,1]$ is $\alpha \times \{-1,1\}$, and
the other $C_j$ remain disjoint from $\alpha \times [-1,1]$.
There is an obvious way to extend this definition of $f$
to a small neighbourhood of $\alpha \times [-1,1]$,
so that it remains unchanged outside of this
neighbourhood. Note that this changes $f$ on 
2-cells $D$ which intersect $\alpha$ only in points,
with the introduction of a new arc of $C_i \cap D$ around
each of these points. (See Figure 1.) Now, we have arranged that
$E_1 - \partial \alpha$ and $E_2 - \partial \alpha$ are both connected.
So, this operation has the effect of
combining $E_1$ and $E_2$ into a single connected
cocycle, thereby reducing the number of components
of $\bigcup C_i$. 

Hence, we may assume that ($\ast$) holds, after
possibly homotoping $f$. This homotopy has the effect
of changing the induced homomorphism
$f_\ast \colon \pi_1(K) \rightarrow \ast^n {\Bbb Z}$
by a conjugacy, but it remains a surjective
homomorphism.

\vskip 18pt
\centerline{\psfig{figure=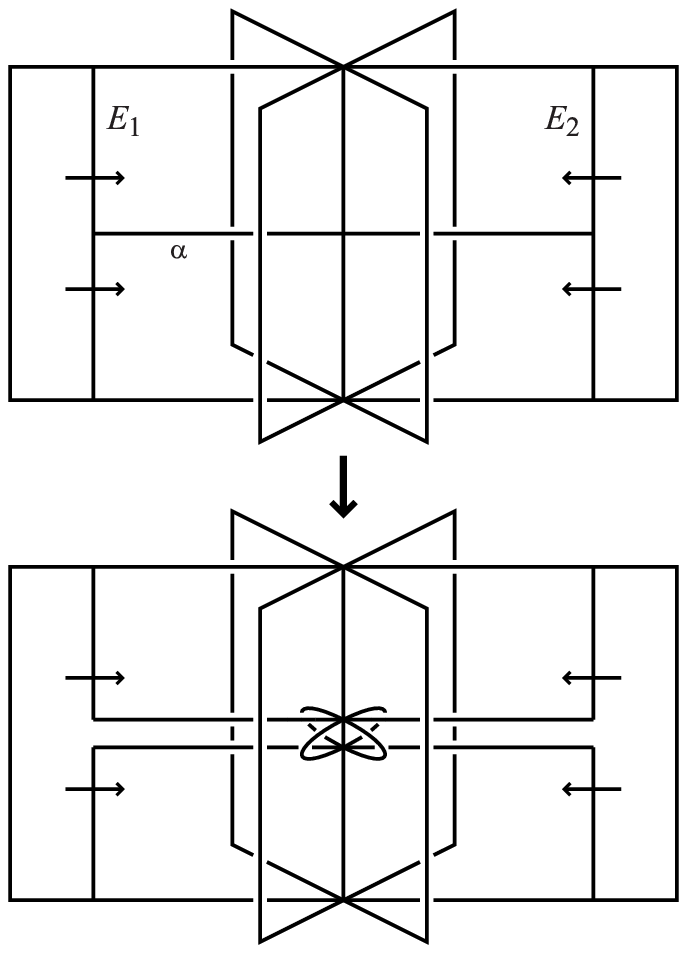,width=2.5in}}
\vskip 24pt
\centerline{Figure 1.}

We claim that $Y$ then has a single vertex, with $n$ edges,
labelled $1, \dots, n$. This will show that
$\bigcup C_i$ is non-separating as required. To prove
this claim, we use the hypothesis that $f_\ast$ is
surjective. This implies that there are loops
$\ell_1, \dots, \ell_n$, based at the basepoint
of $K$, that are sent to the free generators
of $\ast^n {\Bbb Z}$. Pick these loops
so that they have the fewest number of intersections
with $\bigcup C_i$. The loops determine loops
in the graph $Y$. No loop can travel over
$C_i$ in one direction, and then back across
$C_i$ in the other direction. For, by property
($\ast$), it would have to return to the same component
of $C_i$. We could then remove this sub-arc of
the loop, and replace it by an arc in $C_i$, and
then perform a small homotopy, reducing the number
of intersections with $\bigcup C_i$ by two. The
resulting loop still is sent to the same element
of $\ast^n {\Bbb Z}$, which contradicts our
minimality assumption. Hence, the word that $\ell_i$
spells, as it runs over $\bigcup C_i$, is a reduced
word. It therefore runs over $C_i$ exactly once,
and is disjoint from the other cocycles. Hence,
emanating from the vertex of $Y$ that corresponds
to the component of $K - \bigcup C_i$ containing
the basepoint, there is an edge labelled $i$, for
each $i$, and each such edge returns to this vertex. Therefore,
$Y$ is a bouquet of circles, as required. $\square$

In this theorem, we worked with 2-complexes for
convenience. We could just as easily have worked with
smooth manifolds. In this case, transversely oriented,
codimension one submanifolds play the r\^ole of regular 
cocycles. Essentially the same argument as for Theorem 3.1
gives the following.

\noindent {\bf Theorem 3.2.} {\sl Let $M$ be a connected smooth
manifold. Then $\pi_1(M)$ admits a surjective homomorphism
onto $\ast^n {\Bbb Z}$ if and only if $M$ contains $n$
disjoint, transversely oriented, codimension one submanifolds
whose union is non-separating.}

All of the above is fairly well known. What is possibly
less widely known is that one can replicate much of this
work using cohomology with coefficients in ${\Bbb F}_p$,
the field of order a prime $p$. Therefore, fix a prime $p$.

A {\sl regular mod $p$ cocycle} has a similar definition to
a regular cocycle. Again, it is a non-empty finite graph $\Gamma$ 
embedded in $K$, with a little extra structure. It must be disjoint from the
0-skeleton of $K$. However, unlike the case of
regular cocycles, it has two type of vertices,
which we term {\sl edge vertices} and {\sl interior vertices}.
The edge vertices are the intersection of $\Gamma$ with the
1-skeleton of $K$. The vertices of $\Gamma$ on the boundary
of any 2-cell are therefore edge vertices, and we require
them to have valence one in that 2-cell. Each interior vertex
lies in the interior of a 2-cell of $K$.
The edges of $\Gamma$ are
given a transverse orientation and a {\sl weight},
which is a non-zero integer mod $p$. These must satisfy the following local
conditions near the vertices. Near the edge vertices, 
the transverse orientations and the weights must
all be locally equivalent. Around any interior vertex, 
the total weight (signed according to the transverse
orientations) must be congruent to zero mod $p$.
We also insist that each interior vertex has at least
one edge incident to it. (See Figure 2.)

\vskip 18pt
\centerline{\psfig{figure=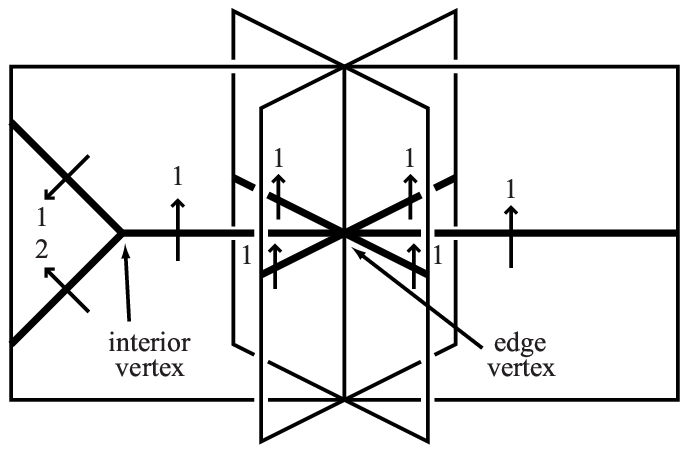}}
\vskip 18pt
\centerline{Figure 2.}

We will see that, as before, any element of $H^1(K; {\Bbb F}_p)$ is
represented by a regular mod $p$ cocycle, and conversely,
a regular mod $p$ cocycle determines a class in 
$H^1(K; {\Bbb F}_p)$. The following states that,
for non-trivial cohomology classes, we may ensure
that the regular mod $p$ cocycle is also non-separating. 

\noindent {\bf Proposition 3.3.} {\sl Let $K$ be a finite
connected 2-complex and let $p$ be a prime. Then any
non-trivial element of $H^1(K; {\Bbb F}_p)$ is represented
by a non-separating regular mod $p$ cocycle.}

\noindent {\sl Proof.} Any element of $H^1(K; {\Bbb F}_p)$
is represented by a cellular 1-dimensional cocycle $c$. This is an assignment
to each oriented 1-cell $e$ of an integer mod $p$ which we
denote by $c(e)$, with
the proviso that the sum of the integers around any
2-cell is zero mod $p$. From this, we build a regular
mod $p$ cocycle $\Gamma$ as follows. Into each 1-cell $e$ for
which $c(e)$ is non-zero mod $p$, we place an
edge vertex of $\Gamma$ with weight $c(e)$. If a 
2-cell contains a 1-cell with non-zero weight in its
boundary, insert into it a single interior vertex. Join
this vertex to each edge vertex in the boundary of
the 2-cell. The fact that the total weight of
$c$ around the 2-cell is zero mod $p$ implies that the
local condition near the interior vertex is satisfied.
Thus, it is trivial
that any element of $H^1(K; {\Bbb F}_p)$ is represented
by a regular mod $p$ cocycle $\Gamma$. 

The aim now is to
ensure that $\Gamma$ is non-separating when the cohomology
class is non-zero. To establish this, we will perform a
sequence of alterations to $\Gamma$. Each will reduce
the number of edge vertices, and so this sequence is guaranteed
to terminate. Suppose that $\Gamma$ is separating,
and let $K_1$ be some component of $K - \Gamma$.
Then, there is some edge vertex in the boundary of $K_1$
that is incident to another component of $K - \Gamma$.
Let $\Gamma'$ be the component of $\Gamma$ minus its
interior vertices that contains this edge vertex.
Then all the edges of $\Gamma'$ are compatibly oriented
and have the same weight $w$, say. Remove $\Gamma'$ from
$\Gamma$. For each edge in the boundary of $K_1$ but not in
$\Gamma'$, add or
subtract $w$ to its weight, according to whether the
transverse orientation of the edge points into or out of
$K_1$. If both sides of the edge lie in $K_1$, then leave
its weight unchanged. If this procedure changes the weight
of any edge to zero mod $p$, then remove it. If any
interior vertices become isolated, remove them. The result
is a new regular mod $p$ cocycle, representing the same
cohomology class, and with fewer edge vertices. Repeating
this process a sufficient number of times, we therefore
end with a non-separating regular mod $p$ cocycle. $\square$

The proof of the above result gives the following extra
information which will be useful later.

\noindent {\bf Addendum 3.4.} {\sl Let $K$ be a finite 2-complex
and let $p$ be a prime. If a regular mod $p$ cocycle
$\Gamma$ represents a non-trivial element of $H^1(K; {\Bbb F}_p)$,
then some subgraph of $\Gamma$ is a regular mod $p$ cocycle (with possibly different weights)
which represents the same cohomology class and is non-separating in $K$.}

There is also a more technical version of Proposition 3.3 that 
deals with subcomplexes.

\noindent {\bf Proposition 3.5.} {\sl Let $K$ be a finite 2-complex
and let $p$ be a prime. Let $L$ be a subcomplex of $K$.
Suppose that there is a non-trivial element $\alpha$ in the kernel
of $H^1(K; {\Bbb F}_p) \rightarrow H^1(L; {\Bbb F}_p)$, the
map induced by inclusion. Then $\alpha$ is represented
by a regular mod $p$ cocycle that is non-separating in $K$
and disjoint from $L$.}

\noindent {\sl Proof.} Pick a cellular cochain $c$ that represents
$\alpha$. Since the restriction of $c$ to $L$ is cohomologically
trivial, it is a coboundary in $L$. Subtracting this coboundary
from $c$ does not change the class it represents, but afterwards
its evaluation on any 1-cell in $L$ is trivial. Thus, when the
construction in the proof of Proposition 3.3 is performed,
a regular mod $p$ cocycle $\Gamma$ is created that is disjoint
from $L$. Applying Addendum 3.4, we can ensure that $\Gamma$
is non-separating in $K$ and still disjoint from $L$. $\square$

There is also a corresponding
version of Theorem 3.1 for regular mod $p$ cocycles, which works best when
$p = 2$. This will be a crucial tool in proving that
certain groups are large.

\noindent {\bf Theorem 3.6.} {\sl Let $K$ be a finite connected 2-complex,
and let $p$ be a prime. Then $\pi_1(K)$
admits a surjective homomorphism onto $\ast^n ({\Bbb Z}/p{\Bbb Z})$
if $K$ contains $n$ disjoint regular mod $p$ cocycles whose union
is non-separating. Furthermore the converse holds
when $p=2$ and $K$ contains no locally separating
points and no 1-cells with valence 1.}

\noindent {\sl Proof.} The proof is very similar to that of
Theorem 3.1, and so we will only focus on those parts where
the details differ.

Suppose first that $K$ contains $n$ disjoint regular mod $p$ cocycles
whose union is non-separating. Then we construct a map
$f \colon K \rightarrow \bigvee^n L(p)$, where $L(p)$ is
the 2-complex consisting of a single 0-cell, a single 1-cell,
and a 2-cell that winds $p$ times around the 1-skeleton.
Outside of a small regular neighbourhood of the cocycles,
everything is sent by $f$ to the central vertex of
$\bigvee^n L(p)$. On product neighbourhoods of the
edges and edge vertices of the cocycles, $f$ is defined to collapse these
products onto an interval and then map this interval $w$ times
around the relevant 1-cell of $\bigvee^n L(p)$, where $w$ is the
weight of the edge.
Finally, near the interior vertices of the cocycles, $f$ maps onto the
relevant 2-cell of $\bigvee^n L(p)$. The proof that
$f_\ast \colon \pi_1(K) \rightarrow \ast^n ({\Bbb Z}/p{\Bbb Z})$
is a surjection is similar to the corresponding
proof for Theorem 3.1.

Suppose now that $\pi_1(K)$ admits a surjective
homomorphism onto $\ast^n ({\Bbb Z}/p{\Bbb Z})$.
Suppose also that $p=2$ and $K$ contains no locally separating
points and no 1-cells with valence 1.
Then, exactly as in the proof of Theorem 3.1,
this homomorphism is induced by a map $f \colon K \rightarrow
\bigvee^n L(2)$. Let $\alpha_i$ be the regular mod $2$ cocycle
in $\bigvee^n L(2)$ that has exactly one edge vertex
in the $i^{\rm th}$ 1-cell and exactly one interior vertex
in the $i^{\rm th}$ 2-cell. Then we may arrange that 
$f^{-1}(\alpha_i)$ forms a regular mod $2$ cocycle $C_i$
for each $i$. We may also arrange that each interior
vertex of $\bigcup C_i$ has valence 2. However, as in the proof of
Theorem 3.1, the union of these cocycles may not yet
be non-separating in $K$. We may need to modify $f$
by a homotopy before this condition is satisfied.

Define a graph $Y$ whose vertices correspond to
complementary components of
$\bigcup C_i$, and whose edges correspond to
the components of $\bigcup C_i$. It may not
be the case that a component of $\bigcup C_i$
has a regular neighbourhood that is a product.
If it is not a product, then using the fact that $p = 2$, it is adjacent
to a single complementary region of $\bigcup C_i$,
and we therefore define the corresponding
edge of $Y$ to be a loop. The
edges of $Y$ come with an integer label between 1
and $n$, depending on which cocycle $C_i$ they
came from. However, they do not necessarily come with 
a well-defined orientation. Again, we will homotope $f$,
to ensure that a certain condition holds:
\item{($\ast'$)} no vertex of $Y$ has two distinct edges adjacent
to it with the same label.

Each modification will reduce the
number of components of $\bigcup C_i$, and
so they are guaranteed to terminate. 
The modifications are exactly as before,
except now the transverse orientations
of $E_1$ and $E_2$ at the endpoints of $\alpha$
might not point towards each other or away from
each other. However, this is easily rectified
by the introduction of two interior vertices
near one of the endpoints of $\alpha$.
The argument now proceeds exactly as in
the proof of Theorem 3.1. $\square$

The following consequence of Theorem 3.6 gives a method for
proving that certain groups are large.

\noindent {\bf Theorem 3.7.} {\sl Let $K$ be a finite cell complex, and let
$A$ and $B$ be subcomplexes such that $K = A \cup B$. Let $p$ be a prime and let ${\Bbb F}_p$
be the field of order $p$. Suppose that both of the maps
$$\eqalign{
H^1(A; {\Bbb F}_p) & \rightarrow H^1(A \cap B; {\Bbb F}_p) \cr
H^1(B; {\Bbb F}_p) & \rightarrow H^1(A \cap B; {\Bbb F}_p) \cr}$$
induced by inclusion are not injections. In the case $p =2$,
suppose also that the kernel of at least one of these maps
has dimension more than one. Then $\pi_1(K)$ admits a
surjective homomorphism onto 
$({\Bbb Z}/p{\Bbb Z}) \ast ({\Bbb Z}/p{\Bbb Z})$.
Furthermore, some normal subgroup of $\pi_1(K)$
with index a power of $p$ admits a surjective
homomorphism onto a non-abelian free group. Hence,
$\pi_1(K)$ is large.}

\noindent {\sl Proof.} We may restrict attention to
the 2-skeleton of $K$, since this has the same fundamental
group as $K$, and since the relevant homomorphisms
between cohomology groups are unchanged. Thus, we
may assume that $K$ is a 2-complex.

Pick a non-trivial element of
the kernel of $H^1(A; {\Bbb F}_p) \rightarrow H^1(A \cap B; {\Bbb F}_p)$.
By Proposition 3.5, this is represented by a regular
mod $p$ cocycle that is disjoint from $A \cap B$ and
that is non-separating in $A$. It is therefore a
regular mod $p$ cocycle in $K$. The same argument
gives a non-separating regular mod $p$ cocycle in
$B$ that is disjoint from $A \cap B$. Hence, we obtain
two disjoint regular mod $p$ cocycles in $K$ whose union is
non-separating. By Theorem 3.6, this implies that
$\pi_1(K)$ admits a surjective homomorphism onto
$({\Bbb Z}/p{\Bbb Z}) \ast ({\Bbb Z}/p{\Bbb Z})$.
When $p$ is odd, 
$({\Bbb Z}/p{\Bbb Z}) \ast ({\Bbb Z}/p{\Bbb Z})$
contains a free non-abelian normal subgroup with index
a power of $p$. The inverse image of this subgroup
in $\pi_1(K)$ is also normal and has index a power of $p$.
It surjects on this free non-abelian group.
This therefore proves the theorem when 
$p$ is odd.

Now, $({\Bbb Z}/2{\Bbb Z})
\ast ({\Bbb Z}/2{\Bbb Z})$ does not have a free non-abelian
group as a subgroup, and so the theorem is not
yet fully proved when $p=2$. In this case, however, we 
are assuming that the kernel of one of the maps,
say $H^1(A; {\Bbb F}_p) \rightarrow H^1(A \cap B ;{\Bbb F}_p)$,
has dimension at least two. Construct the finite-sheeted
covering space of $A$ corresponding to this kernel.
The inverse image of $A \cap B$ is a disjoint union of 
at least four copies of $A \cap B$. Attach to each of
these a copy of $B$. The result is a finite-sheeted
regular cover $\tilde K$ of $K$ with degree a power of $2$. 
In each copy of $B$, there 
is a non-separating regular mod $2$ cocycle.
The union of these is therefore non-separating in $\tilde K$. 
Thus, by Theorem 3.6, $\pi_1(\tilde K)$ admits a surjective homomorphism
onto $\ast^4 ({\Bbb Z}/2{\Bbb Z})$. This contains a
normal free non-abelian subgroup, with index a power of $2$.
Its inverse image in $\pi_1(\tilde K)$ surjects onto
this non-abelian free group. By passing a further subgroup
if necessary, we may assume that this is normal in $\pi_1(K)$
and has index a power of $2$ in $\pi_1(K)$. $\square$

Thus, one route to proving that a cell complex $K$
has large fundamental group is to find a decomposition
into subcomplexes $A$ and $B$ where $|H_1(A; {\Bbb F}_p)|$ and 
$|H_1(B; {\Bbb F}_p)|$
are both bigger than $2|H_1(A \cap B; {\Bbb F}_p)|$. This suggests the
following definition.

\noindent {\bf Definition.} Let $K$ be a finite cell complex.
Consider all ways of decomposing $K$ into two sets
$A$ and $B$, where $A$ and $B$ are subcomplexes in
some subdivision of the cell structure on $K$. Let
the {\sl mod $p$ Cheeger constant} of $K$, denoted
$h_p(K)$, be
$$\inf \left \{ {|H_1(A \cap B; {\Bbb F}_p)| \over
\min \{ |H_1(A; {\Bbb F}_p)|, |H_1(B; {\Bbb F}_p)| \} }
\right \}.$$

Theorem 3.7 has the following immediate corollary.

\noindent {\bf Corollary 3.8.} {\sl Let $K$ be
a finite connected  cell complex, and let $p$ be a prime.
Suppose that
$$h_p(K) < \cases{ 1 & if $p$ is odd; \cr
1/2 & if $p = 2$.}$$
Then $\pi_1(K)$ admits a
surjective homomorphism onto 
$({\Bbb Z}/p{\Bbb Z}) \ast ({\Bbb Z}/p{\Bbb Z})$.
Furthermore, some normal subgroup of $\pi_1(K)$
with index a power of $p$ admits a surjective
homomorphism onto a non-abelian free group. Hence,
$\pi_1(K)$ is large.}

The following result summarises much of what has been done
in this section.

\noindent {\bf Theorem 3.9.} {\sl Let $K$ be a finite connected 2-complex
with fundamental group $G$. Suppose that $K$ has no
locally separating points and no 1-cells with valence 1.
Then the following are equivalent:
\item{(i)} $G$ is large;
\item{(ii)} in some finite-sheeted covering space $\tilde K$ of
$K$, there are two disjoint regular cocycles whose union is non-separating;
\item{(iii)} in some finite-sheeted covering space $\tilde K$ of
$K$, there are two disjoint regular mod $p$ cocycles 
whose union is non-separating, for some odd prime $p$;
\item{(iv)} in some finite-sheeted covering space $\tilde K$ of
$K$, there are three disjoint regular mod $2$ cocycles 
whose union is non-separating.

}

\noindent {\sl Proof.} Note first that condition of
having no locally separating points and no 1-cells
with valence 1 is preserved under finite covers.

(i) $\Rightarrow$ (ii): Since $G$ is
large, some finite index subgroup of $G$ admits a surjective
homomorphism onto ${\Bbb Z} \ast {\Bbb Z}$. Let $\tilde K$
be the covering space of $K$ corresponding to this subgroup.
By Theorem 3.1, it has two disjoint regular cocycles
whose union is non-separating.

(ii) $\Rightarrow$ (iii):
This is obvious, because a regular cocycle becomes a regular
mod $p$ cocycle when every edge is given weight $1$.

(iii) $\Rightarrow$ (i): By Theorem 3.6, the fundamental
group of $\tilde K$ admits a surjective homomorphism onto
$({\Bbb Z}/p{\Bbb Z}) \ast ({\Bbb Z}/p{\Bbb Z})$. But this
contains a non-abelian free group as a finite index normal subgroup.
Hence, $G$ is large.

(i) $\Rightarrow$ (iv): This is very similar to (i) $\Rightarrow$ (ii)
$\Rightarrow$ (iii). Since $G$ is large, some finite index subgroup
admits a surjective homomorphism onto ${\Bbb Z} \ast {\Bbb Z} \ast {\Bbb Z}$.
The corresponding covering space has three disjoint
regular cocycles whose union is non-separating. Each is,
by definition, a regular mod $2$ cocycle when every edge is
given weight 1.

(iv) $\Rightarrow$ (i): This proof is essentially the same
as (iii) $\Rightarrow$ (i),
using the fact that $({\Bbb Z}/2{\Bbb Z}) \ast ({\Bbb Z}/2{\Bbb Z})
\ast ({\Bbb Z}/2{\Bbb Z})$ has a free non-abelian group as
a finite index normal subgroup. $\square$

\vskip 18pt
\centerline{\caps 4. Cheeger decompositions of coset diagrams}
\vskip 6pt

The following result was a key technical lemma in [3]
(Lemma 2.1 there).

\noindent {\bf Lemma 4.1.} {\sl Let $X$ be a Cayley graph
of a finite group, and let $D$ be a non-empty subset
of $V(X)$ such that $|\partial D|/|D| = h(X)$ and
$|D| \leq |V(X)|/2$. Then
$|D| > |V(X)|/4$. Furthermore, the subgraphs induced
by $D$ and its complement $D^c$ are connected.}

This was useful when analysing finite index normal subgroups $H$
of a group $G$, because then a finite generating set for $G$ determines
a Cayley graph of $G/H$. However, in this paper, we wish to consider
subgroups that are not necessarily normal.
Thus, the following generalisation will be necessary.

\noindent {\bf Proposition 4.2.} {\sl Let $G$ be a group with
a finite generating set $S$, and let $\{ G_i \}$ be
a sequence of finite index subgroups, where each $G_i$ is normal in
$G_{i-1}$. Let $S$ be a finite set of generators
for $G$, and let $X_i$ be $X(G/G_i;S)$. Then $h(X_i)$ is a non-increasing
sequence. Suppose that, for some $i$, $h(X_i) < h(X_{i-1})$.
Then, there is some non-empty subset $D$
of $V(X_i)$ such that $|\partial D|/|D| = h(X_i)$ and
$|V(X_i)|/4 < |D| \leq |V(X_i)|/2$.}

\noindent {\sl Proof.} The fact that $h(X_i)$ is non-increasing
is trivial. Therefore, let us concentrate on the second
part of the proposition. Consider a non-empty subset $D$ 
of $V(X_i)$ such that
$|\partial D|/|D| = h(X_i)$ and $|D| \leq |V(X_i)|/2$.
Pick $D$ so that $|D|$ is as large as possible subject
to these two conditions. Let us suppose
that $|D| \leq |V(X_i)|/4$, with the aim
of reaching a contradiction. Now, $G_i$ is normal in $G_{i-1}$
and so $G_{i-1}/G_i$ acts on $X_i$ by covering transformations.
Let $g$ be any element of $G_{i-1}/G_i$. We consider $g(D) \cup D$.
It is shown in [3] (see the proof of Lemma 2.1 there) that
$$|\partial(g(D) \cup D)| = |\partial D| + |\partial g(D)| - 
|\partial (g(D) \cap D)| - 2 e(g(D) - D,D - g(D)),$$
where $e(g(D) - D, D - g(D))$ denotes the number of edges joining
$g(D) - D$ to $D - g(D)$. By the definition of
$h(X_i)$, we must have that $|\partial (g(D) \cap D)| \geq h(X_i) |g(D) \cap D|$. Thus,
$$|\partial (g(D) \cup D)| \leq h(X_i) (|D| + |g(D)| - |g(D) \cap D|) = h(X_i) |g(D) \cup D|.$$
Now, $g(D) \cup D$ is at most half the vertices of $X_i$,
by our assumption that $|D| \leq |V(X_i)|/4$. As $|D|$
was assumed to be maximal, $|g(D) \cup D|$ must be equal to $|D|$ and hence $g(D) = D$. This is true for each $g \in G_{i-1}/G_i$.
Thus, $D$ is invariant under the action of $G_{i-1}/G_i$ on $X_i$,
and therefore descends to a subset $D'$ of $V(X_{i-1})$. 
Now, $|\partial D'| = |\partial D|/[G_{i-1}:G_i]$ and $|D'| = |D|/[G_{i-1}:G_i]$.
Hence, 
$$h(X_{i-1}) \leq |\partial D'|/|D'| = |\partial D|/|D| = h(X_i) \leq h(X_{i-1}).$$
Thus, these must be equalities, which contradicts our
hypothesis that $h(X_i) < h(X_{i-1})$. Hence, it must have
been the case that $|D| > |V(X_i)|/4$. $\square$

\vskip 18pt
\centerline{\caps 5. Proof of the main theorem }
\vskip 6pt

In this paper, we will be concentrating on groups $G$ having
a sequence of finite index subgroups $\{G_i \}$ with
linear growth of mod $p$ homology, for some prime $p$.
It will be helpful to introduce a quantity that measures
the growth rate of $d_p(G_i)$.
This is the {\sl mod $p$ homology gradient} which is
defined to be 
$$\inf_i {(d_p(G_i) -1) \over [G:G_i]}.$$
This quantity is most relevant when each $G_{i+1}$ is
normal in $G_i$ and has index a power
of $p$. In this case, we have the following well-known proposition.

\noindent {\bf Proposition 5.1.} {\sl Let $G$ be a finitely
generated group, and let $H$ be a subnormal subgroup with index
a power of a prime $p$. Then
$$d_p(H) - 1 \leq [G:H] (d_p(G) - 1).$$}

This appears as Proposition 3.7 in [5]
for example. It implies that when each $G_{i+1}$ is
normal in $G_i$ and has index a power
of $p$, $(d_p(G_i) - 1)/[G:G_i]$
is a non-increasing function of $i$. In particular, the
infimum in the definition of mod $p$ homology gradient is
a limit.

We will, in fact, need the following stronger result.

\noindent {\bf Proposition 5.2.} {\sl Let $K$ be a connected 
2-complex, and let $\Gamma$ be a connected union of 1-cells
such that the map $H_1(\Gamma; {\Bbb F}_p)
\rightarrow H_1(K; {\Bbb F}_p)$ induced by inclusion is 
a surjection, for some prime
$p$. Let $\tilde K \rightarrow K$ be a finite-sheeted
covering such that $\pi_1(\tilde K)$ is subnormal
in $\pi_1(K)$ and has index a power of $p$.
Let $\tilde \Gamma$ be the inverse image of
$\Gamma$ in $\tilde K$. Then $\tilde \Gamma$ is
connected and the map $H_1(\tilde \Gamma; {\Bbb F}_p)
\rightarrow H_1(\tilde K; {\Bbb F}_p)$ induced by inclusion
is a surjection.}

To prove this, we will require the following.

\noindent {\bf Lemma 5.3.} {\sl Let $\Gamma$ be a path-connected
subset of a path-connected space $L$ such that the map $\pi_1(\Gamma) \rightarrow 
\pi_1(L)$ induced by inclusion is surjection. Let
$q \colon \tilde L \rightarrow L$ be a covering map,
and $\tilde \Gamma$ be the inverse image of $\Gamma$
in $\tilde L$. Then $\tilde \Gamma$ is
path-connected and the map $\pi_1(\tilde \Gamma)
\rightarrow \pi_1(\tilde L)$ induced by inclusion
is a surjection.}

\noindent {\sl Proof.} Let $b$ be a basepoint for
$L$ in $\Gamma$. The restriction of $q$ to any path-component
of $\tilde \Gamma$ is a covering map onto $\Gamma$,
which is therefore surjective. Thus,
to show that $\tilde \Gamma$
is path-connected, it suffices to show that any two
points of $q^{-1}(b)$ lie in the same path-component of
$\tilde \Gamma$. We may assume that one of these
points is a basepoint $\tilde b$ of $\tilde L$.
Pick a path from $\tilde b$ to the other point
in $q^{-1}(b)$. This projects to a loop $\ell $ in 
$L$ based at $b$. Since $\pi_1(\Gamma,b) \rightarrow \pi_1(L,b)$
is a surjection, $\ell$ is homotopic, relative
to its endpoints, to a loop in $\Gamma$. This lifts
to a path in $\tilde \Gamma$ joining the two
points of $q^{-1}(b)$.

We now show that $\pi_1(\tilde \Gamma, \tilde b)
\rightarrow \pi_1(\tilde L, \tilde b)$ is
a surjection. Given any loop $\tilde \ell$
in $\tilde L$ based at $\tilde b$, we project
it to a loop $\ell$ in $L$. This is homotopic
relative to its endpoints to a loop in $\Gamma$.
This homotopy lifts to a homotopy, relative
to endpoints, between $\tilde \ell$ and a loop in
$\tilde \Gamma$. $\square$

\noindent {\sl Proof of Proposition 5.2.} Note first that an
obvious induction allows us to reduce to the
case where $\pi_1(\tilde K)$ is a normal
subgroup of $\pi_1(K)$ with index a power of $p$.

Pick a maximal tree in $\Gamma$
and extend it to a maximal tree $T$ in the 1-skeleton
of $K$. Let $\overline K$ be obtained from $K$ by collapsing
$T$ to a point, and let $\overline \Gamma$ be the
image of $\Gamma$ in $\overline K$. Then clearly
the map $H_1(\overline \Gamma; {\Bbb F}_p) \rightarrow
H_1(\overline K; {\Bbb F}_p)$ induced by inclusion
is a surjection. Suppose that we could prove the
theorem for $\overline K$ and $\overline \Gamma$.
Then this would clearly imply the theorem for
$K$ and $\Gamma$. Thus, we may assume that
$K$ has a single 0-cell. It therefore specifies
a presentation for $\pi_1(K)$, once we have picked
an orientation on each of the 1-cells of $K$.

Let $G$ and $H$ denote the groups $\pi_1(K)$
and $\pi_1(\tilde K)$ respectively.
Let $H'$ denote $[H,H]H^p$,
the subgroup of $H$ generated by the commutators and $p^{\rm th}$ powers
of $H$. This is a
characteristic subgroup of $H$, with index a power
of $p$. We are assuming that $H$ is a normal subgroup
of $G$ with index a power of $p$. Hence, $H'$ is
a normal subgroup of $G$ with index a power of $p$.
In other words, $G/H'$ is a finite $p$-group.

Now, $H_1(G/H'; {\Bbb F}_p)$ is isomorphic to $H_1(G; {\Bbb F}_p)$. 
Hence, the 1-cells of $\Gamma$
form a generating set for $H_1(G/H'; {\Bbb F}_p)$.
It is a well known fact that in any finite $p$-group $C$,
a set of elements forms a generating set for $C$
if and only if they form a generating set for
$H_1(C; {\Bbb F}_p)$. Thus, the 1-cells of $\Gamma$ form
a generating set for $G/H'$. Let $L$ be the 2-complex obtained
from $K$ by attaching a 2-cell along each word in
$H'$. Then $L$ has fundamental group $G/H'$. The map
$\pi_1(\Gamma) \rightarrow \pi_1(L)$ induced by inclusion
is a surjection.
Let $\tilde L$ be the covering space of $L$ corresponding to
the subgroup $H/H'$. This is obtained from $\tilde K$
by attaching various 2-cells. But one may view
their 1-skeletons as the same. By Lemma 5.3, the inverse image 
of $\Gamma$ in $\tilde L$
is a connected graph. This is a copy of $\tilde \Gamma$,
and so $\tilde \Gamma$ is connected.
The map $\pi_1(\tilde \Gamma)
\rightarrow \pi_1(\tilde L)$ induced by inclusion is
a surjection, by Lemma 5.3. The natural map $\pi_1(\tilde L)
\rightarrow H_1(\tilde L; {\Bbb F}_p)$ is a surjection.
This implies that the map $H_1(\tilde \Gamma; {\Bbb F}_p) \rightarrow
H_1(\tilde L; {\Bbb F}_p)$ is a surjection.
The map $H_1(\tilde K; {\Bbb F}_p) \rightarrow H_1(\tilde L;
{\Bbb F}_p)$ induced by inclusion is an isomorphism.
Hence, $H_1(\tilde \Gamma; {\Bbb F}_p) \rightarrow
H_1(\tilde K; {\Bbb F}_p)$ is a surjection, as required.
$\square$

Before we prove Theorem 1.1, we introduce some terminology.
If $K$ is a topological space and $p$ is a prime, then
$d_p(K)$ denotes the dimension of $H_1(K; {\Bbb F}_p)$.

\noindent {\sl Proof of Theorem 1.1.} Suppose that
$\{ G_i \}$ has linear growth of mod $p$ homology,
and that $G$ does not have Property $(\tau)$ with
respect to $\{G_i\}$. Our aim is to show that
some $G_i$ admits a surjective homomorphism onto
$({\Bbb Z}/p{\Bbb Z}) \ast ({\Bbb Z}/p{\Bbb Z})$ and 
that some normal subgroup of $G_i$, with index
a power of $p$, admits a surjective homomorphism onto
a non-abelian free group.

We fix $\epsilon$ to be some real number strictly
between $0$ and $\sqrt{10}/3 -1$, but where we
view it as very small.
Since the mod $p$ homology gradient of $\{ G_i \}$ is non-zero,
there is some $j$ such that $(d_p(G_j)-1)/[G:G_j]$ 
is at most $(1 + \epsilon)$ times the mod $p$ homology gradient of $\{ G_i \}$.
The mod $p$ homology gradient of $\{ G_i : i \geq j \}$ (viewed
as subgroups of $G_j$) is $[G:G_j]$ times the mod $p$ homology
gradient of $\{G_i\}$ (viewed as subgroups of $G$).
So, $d_p(G_j) - 1$ is at most $(1 + \epsilon)$ times the
mod $p$ homology gradient of $\{ G_i : i \geq j \}$.
Hence, by replacing $G$ by $G_j$,
and replacing $\{ G_i \}$ by $\{ G_i : i \geq j \}$,
we may assume that $d_p(G)-1$ is at most $(1 + \epsilon)$ times the
mod $p$ homology gradient of $\{ G_i \}$. We may also
assume (by replacing $G$ by $G_1$) that the index of
each $G_i$ in $G$ is a power of $p$.

Let $S$ be a set of elements of $G$ that forms
a basis for $H_1(G; {\Bbb F}_p)$.
Extend this to a finite generating set $S_+$ for $G$.
Let $K$ be a finite 2-complex having 
fundamental group $G$, arising from a finite
presentation of $G$ with generating set $S_+$. 
Thus, $K$ has a single vertex and $|S_+|$ edges. Let $L$ be
the sum of the lengths of the relations in this
presentation.
Let $K_i \rightarrow K$ be the covering 
corresponding to $G_i$. Our aim is to show
that its mod $p$ Cheeger constant satisfies
the inequality $h_p(K_i) < 1/2$ for all sufficiently
large $i$. Corollary 3.8 will then prove the 
theorem.

Let $X_i$ be the
1-skeleton of $K_i$. Then $X_i = X(G/G_i; S_+)$.
Let $\Gamma_i$ be the subgraph
of $X_i$ consisting of those edges labelled by $S$.
By Proposition 5.2, $\Gamma_i$ is connected and
the inclusion $\Gamma_i \rightarrow K_i$ induces
a surjection $H_1(\Gamma_i; {\Bbb F}_p) \rightarrow H_1(K_i; {\Bbb F}_p)$.

Since we are assuming that
$G$ does not have Property $(\tau)$ with respect
to $\{ G_i \}$, $\inf_i h(X_i) = 0$.
Since the subgroups $G_i$ are nested, $h(X_i)$ is
a non-increasing sequence. Hence $h(X_i) \rightarrow 0$.
Let us focus on those values of $i$ for which
$h(X_i) < h(X_{i-1})$. This occurs infinitely often. 
Proposition 4.2 asserts that there is a non-empty subset $D_i$
of $V(X_i)$ such that $|\partial D_i|/|D_i| = h(X_i)$ and
$|V(X_i)|/4 < |D_i| \leq |V(X_i)|/2$. We will use $D_i$
to construct a decomposition of $K_i$ into two
overlapping subsets. Let $A_i$ (respectively, $B_i$) 
be the closure of the union of those cells in $K_i$ 
that intersect $D_i$ (respectively, $D_i^c$).
Let $C_i$ be $A_i \cap B_i$.
The edges of $A_i \cap \Gamma_i$ are of three
types (that are not mutually exclusive):
\item{(i)} those edges with both endpoints in $D_i$,
\item{(ii)} those edges in $\partial D_i$,
\item{(iii)} those edges in the boundary of a 2-cell that
intersects both $D_i$ and $D_i^c$.

\noindent If we consider the $d_p(G)$ oriented edges of $\Gamma_i$
emanating from each vertex in $D_i$,
we will cover every edge in (i),
and possibly others. Hence, there are at most
$|D_i| d_p(G)$ edges of type (i) in $A_i \cap \Gamma_i$.

Any type (iii) edge lies in a 2-cell that
intersects both $D_i$ and $D_i^c$. This 2-cell therefore intersects
an edge in $\partial D_i$. Consider one of the endpoints of
the latter edge. At most $L$ 2-cells run over this vertex. Each 2-cell
runs over at most $L$ edges. So, there are no more 
than  $|\partial D_i|L^2$
type (iii) edges. There are $|\partial D_i|$
type (ii) edges, and so, there are
at most $|\partial D_i| (L^2 + 1)$ type (ii) and (iii) edges
in total. A similar argument gives that there
are at most $|\partial D_i|(L^2 + 2)$ vertices in $C_i$.

We claim that each component of $A_i \cap \Gamma_i$ and $B_i \cap \Gamma_i$
contains a vertex in $C_i$. Consider any component of $A_i \cap \Gamma_i$. Since
$\Gamma_i$ is connected, there is a path in $\Gamma_i$
from this component to $B_i \cap \Gamma_i$. The first point in
this path that lies in $B_i$ is the required vertex in $C_i$.
The argument for components of $B_i \cap \Gamma_i$ is similar.
So, $|A_i \cap \Gamma_i|$ and $|B_i \cap \Gamma_i|$ are both
at most $|\partial D_i|(L^2 + 2)$.

Now, the following is an excerpt from the Mayer-Vietoris sequence applied to 
$\Gamma_i \cap A_i$ and $\Gamma_i \cap B_i$:
$$H_1(\Gamma_i \cap A_i; {\Bbb F}_p) \oplus H_1(\Gamma_i \cap B_i ; {\Bbb F}_p)
\rightarrow H_1(\Gamma_i; {\Bbb F}_p) \rightarrow H_0(\Gamma_i \cap C_i; {\Bbb F}_p).$$
The exactness of this sequence implies
that the subspace of $H_1(\Gamma_i; {\Bbb F}_p)$ 
generated by the images of $H_1(\Gamma_i \cap A_i; {\Bbb F}_p)$
and $H_1(\Gamma_i \cap B_i ; {\Bbb F}_p)$ has codimension
at most the number of components of $\Gamma_i \cap C_i$.
This is at most the number of vertices in $C_i$,
which is at most $|\partial D_i|(L^2+2) $. Let 
${\rm Im}(H_1(\Gamma_i \cap A_i; {\Bbb F}_p))$ denote the
image of $H_1(\Gamma_i \cap A_i; {\Bbb F}_p)$ in $H_1(K_i; {\Bbb F}_p)$,
and define ${\rm Im}(H_1(\Gamma_i \cap B_i; {\Bbb F}_p))$
and ${\rm Im}(H_1(\Gamma_i; {\Bbb F}_p))$ similarly.
Note that this latter group is all of
$H_1(K_i; {\Bbb F}_p)$ by Proposition 5.2.
We deduce that the sum of the subspaces
${\rm Im}(H_1(\Gamma_i \cap A_i; {\Bbb F}_p))$
and 
${\rm Im}(H_1(\Gamma_i \cap B_i; {\Bbb F}_p))$
has codimension at most $|\partial D_i|(L^2+2)$ in $H_1(K_i ; {\Bbb F}_p)$.

\vfill\eject
Now, $\Gamma_i \cap A_i$ has at most 
$|D_i| d_p(G) + |\partial D_i| (L^2 + 1)$ edges.
It has at least $|D_i|$ vertices. Hence,
$$\eqalign{
d_p(\Gamma_i \cap A_i) 
&= -\chi(\Gamma_i \cap A_i) + |\Gamma_i \cap A_i| \cr
&\leq |D_i| d_p(G) + |\partial D_i| (L^2 + 1) - |D_i| + |\partial D_i|(L^2 + 2) \cr
&= |D_i| (d_p(G) -1 + h(X_i) (2L^2+3)) \cr
&\leq \textstyle{1 \over 2} [G:G_i] (d_p(G) -1 + h(X_i) (2L^2+3)) \cr
&\leq \textstyle{1 \over 2} (1 + \epsilon) [G:G_i] (d_p(G) -1) 
\hbox{ when }h(X_i) \hbox{ is sufficiently small} \cr
&\leq \textstyle{1 \over 2} (1 + \epsilon)^2 (d_p(G_i) - 1).}
$$
A similar sequence of inequalities holds for
$d_p(\Gamma_i \cap B_i)$ but with $|D_i|$ replaced
throughout by $|D_i^c|$ and with ${1 \over 2}$
replaced throughout by ${3 \over 4}$. Here,
we are using the fact that
$|D_i^c| \leq {3 \over 4} [G:G_i]$. So, when $h(X_i)$ is sufficiently
small, ${\rm Im}(H_1(\Gamma_i \cap A_i; {\Bbb F}_p))$
and ${\rm Im}(H_1(\Gamma_i \cap B_i; {\Bbb F}_p))$ each
have dimension at most ${3 \over 4} (1 + \epsilon)^2 d_p(G_i)$.
Note that ${3 \over 4}(1 + \epsilon)^2 < {5 \over 6}$,
by our assumption 
that $\epsilon < \sqrt{10}/3 -1$. We saw
above that the sum of ${\rm Im}(H_1(\Gamma_i \cap A_i; {\Bbb F}_p))$
and ${\rm Im}(H_1(\Gamma_i \cap B_i; {\Bbb F}_p))$ has codimension
at most $|\partial D_i|(L^2+2)$, which equals
$h(X_i) |D_i| (L^2 + 2)$, and this is small compared
with $d_p(G_i)$. Therefore, when $h(X_i)$ is sufficiently
small, ${\rm Im}(H_1(\Gamma_i \cap A_i; {\Bbb F}_p))$
and ${\rm Im}(H_1(\Gamma_i \cap B_i; {\Bbb F}_p))$
each have dimension at least $d_p(G_i)/6$.
Since $H_1(\Gamma_i \cap A_i; {\Bbb F}_p) \rightarrow
H_1(K_i; {\Bbb F}_p)$ factors through
$H_1(A_i; {\Bbb F}_p)$, this must also have
dimension at least $d_p(G_i)/6$. 
When $h(X_i)$ is sufficiently
small, this is significantly more than $d_p(C_i)$.
Thus, we deduce that, when $i$ is sufficiently
large, $d_p(C_i)$ is less than both $d_p(A_i) - 1$
and $d_p(B_i) - 1$. The mod $p$ Cheeger constant
of $K_i$ is therefore less than $1/2$. Corollary 3.8
then implies that $G_i$ admits a
surjective homomorphism onto 
$({\Bbb Z}/p{\Bbb Z}) \ast ({\Bbb Z}/p{\Bbb Z})$.
Furthermore, some normal subgroup of $G_i$
with index a power of $p$ admits a surjective
homomorphism onto a non-abelian free group. Hence,
$G$ is large. $\square$

\vskip 18pt
\centerline{\caps 6. Error-correcting codes and large groups}
\vskip 6pt

Let $G$ be a finitely presented group, and let $\{ G_i \}$ be a 
nested sequence of finite index subgroups. Suppose that $\{ G_i \}$ has
linear growth of mod $p$ homology. Does this imply
that $G$ is large? Let $K$ be a
finite 2-complex with fundamental group $G$, and let $K_i$
be the covering space corresponding to the subgroup $G_i$.
Then one might suspect that the sheer number of elements
of $H^1(K_i; {\Bbb F}_p)$ might force the existence of
two regular mod $p$ cocycles that are disjoint and whose union is
non-separating. Hence, by Theorem 3.6,
$G_i$ would admit a surjective homomorphism onto
$({\Bbb Z}/p{\Bbb Z}) \ast ({\Bbb Z}/p{\Bbb Z})$,
establishing (i), at least when $p$ is odd. However, it appears not to be possible
to turn this reasoning into a proof, due to the
intervention of error-correcting codes. In this section, we explain
how these codes play a r\^ole.

We first introduce a new concept: the {\sl relative size} of
a cohomology class. Let $K$ be a finite cell complex. For a cellular
1-dimensional cocycle $c$ on $K$, let its support ${\rm supp}(c)$
be those 1-cells with non-zero evaluation under $c$.
For an element $\alpha \in H^1(K; {\Bbb F}_p)$,
consider the following quantity. The {\sl relative size} of $\alpha$ is
$${\min \{ |{\rm supp}(c)|: c \hbox{ is a cellular cocycle 
representing } \alpha \} \over 
\hbox{Number of 1-cells of } K}.$$
The relevance of this quantity is apparent in the following
result.

\noindent {\bf Theorem 6.1.} {\sl Let $K$ be a finite connected
2-complex, and let $\{ K_i \rightarrow K\}$ be a
collection of finite-sheeted covering spaces.
Suppose that $\{ \pi_1(K_i) \}$ has linear
growth of mod $p$ homology for some prime $p$.
Then one of the following must hold:
\item{(i)} $\pi_1(K_i)$ admits a surjective
homomorphism onto $({\Bbb Z}/p{\Bbb Z}) \ast ({\Bbb Z}/p{\Bbb Z})$
for infinitely many $i$, and $\pi_1(K)$ is large, or
\item{(ii)} there is some $\epsilon > 0$ such that
the relative size of any non-trivial class in $H^1(K_i; {\Bbb F}_p)$
is at least $\epsilon$, for all $i$.

}

The following will be useful in the proof of this.

\noindent {\bf Lemma 6.2.} {\sl Let $K$ be a finite
2-complex. Let $M$ be the maximal valence of any
its 1-cells. Let $c$ be a cellular cocycle representing
a class $\alpha$ in $H^1(K; {\Bbb F}_p)$, for some
prime $p$. Then $\alpha$ is represented by a regular
mod $p$ cocycle $\Gamma$ containing at most
$M|{\rm supp}(c)|$ edges, at most $|{\rm supp}(c)|$ edge vertices.}

\noindent {\sl Proof.} Recall from Proposition 3.3 the 
construction of a regular mod $p$ cocycle $\Gamma$ from
the cellular cocycle $c$. Each 1-cell in
${\rm supp}(c)$ is assigned an
edge vertex of $\Gamma$. Each such edge vertex is
adjacent to at most $M$ edges. Also, every
edge is adjacent to some edge vertex. 
Hence, $\Gamma$ contains  at most
$M|{\rm supp}(c)|$ edges. $\square$

\noindent {\sl Proof of Theorem 6.1.} 
Suppose that (ii) does not hold. Then there exist 
non-trivial elements of $H^1(K_i; {\Bbb F}_p)$ 
with arbitrarily small relative size. Let $\Gamma$ be a
regular mod $p$ cocycle representing one of these
cohomology classes, and let $e(\Gamma)$ and $ev(\Gamma)$ denote its
number of edges and edge vertices. By Lemma 6.2, we may ensure
that the ratios of 
$e(\Gamma)$ and $ev(\Gamma)$ to the number of 1-cells of $K_i$
is arbitrarily close to zero. 
(Note that the maximal valence of the 1-cells of $K_i$ is the
same for all $i$.)
By Addendum 3.4, we may arrange that $\Gamma$ is non-separating,
without increasing its number of edges and edge vertices.
Note that $e(\Gamma)$ forms an upper bound on $d_p(\Gamma)$.
Let $N(\Gamma)$ be a thin regular neighbourhood
of $\Gamma$. Then $\partial N(\Gamma)$ is a graph with
as many edges as $\Gamma$, and at most $2 ev(\Gamma)$
components. Thus, $d_p(\partial N(\Gamma))$
is bounded above by $e(\Gamma)$.
We are assuming that $\{ \pi_1(K_i) \}$ has
linear growth of mod $p$ homology. Hence, the ratios
of $e(\Gamma)$ and $ev(\Gamma)$ to $d_p(K_i)$ are both arbitrarily
close to zero.  Consider the Mayer-Vietoris sequence
applied to $N(\Gamma)$ and $K_i - {\rm int}(N(\Gamma))$:
$$\eqalign{
H_1(\partial N(\Gamma); {\Bbb F}_p) &\rightarrow 
H_1(N(\Gamma); {\Bbb F}_p) \oplus H_1(K_i - {\rm int}(N(\Gamma)); {\Bbb F}_p) \cr
&\qquad \rightarrow H_1(K_i; {\Bbb F}_p) \rightarrow H_0(\partial N(\Gamma); {\Bbb F}_p).}$$
Now, the dimensions of $H_1(\partial N(\Gamma); {\Bbb F}_p)$,
$H_1(N(\Gamma); {\Bbb F}_p)$
and $H_0(\partial N(\Gamma); {\Bbb F}_p)$ are all small
compared with $d_p(K_i)$. Hence, the ratio of
$d_p(K_i)$ and $d_p(K_i - {\rm int}(N(\Gamma)))$ tends
to 1. So, the ratio of $d_p(\partial N(\Gamma))$
and $d_p(K_i - {\rm int}(N(\Gamma)))$ tends to zero.
Therefore, the map $H^1(K_i - {\rm int}(N(\Gamma)); {\Bbb F}_p) \rightarrow 
H^1(\partial N(\Gamma); {\Bbb F}_p)$ induced by inclusion
has non-trivial kernel. Subdivide $K_i$ so that $N(\Gamma)$
is a subcomplex. By Proposition 3.5, there
is a regular mod $p$ cocycle $\Gamma'$ in 
$K_i - {\rm int}(N(\Gamma))$ such that $\Gamma'$ is non-separating
in $K_i - {\rm int}(N(\Gamma))$ and disjoint from $\partial N(\Gamma)$. So,
$\Gamma \cup \Gamma'$ is non-separating in $K_i$.
By Theorem 3.6, $\pi_1(K_i)$
admits a surjective homomorphism onto 
$({\Bbb Z}/p{\Bbb Z}) \ast ({\Bbb Z}/p{\Bbb Z})$.
When $p >2$, this gives (i). So, let us suppose now
that $p = 2$. We may assume that the kernel of
$H^1(K_i - {\rm int}(N(\Gamma)); {\Bbb F}_p) \rightarrow 
H^1(\partial N(\Gamma); {\Bbb F}_p)$ has dimension at least two.
Pick two linearly independent elements in this
kernel, and consider the cover of $K_i - {\rm int}(N(\Gamma))$,
with order $4$, dual to these two elements. This extends to a
cover $\tilde K_i$ of $K_i$.
The inverse image of $\Gamma$ in $\tilde K_i$ has
at least 4 components. The complement of their union is, by construction, connected.
So, by Theorem 3.6, $\pi_1(\tilde K_i)$ admits a surjective
homomorphism onto $\ast^4 ({\Bbb Z}/2{\Bbb Z})$,
and hence $G$ is large. $\square$

Theorem 6.1 leads naturally to the following question:
how can the relative sizes of the non-trivial cohomology classes 
of $K_i$ not have zero infimum?
The answer is: when they form error correcting codes
with large Hamming distance.

Recall that a {\sl linear code} is a subspace ${\cal C}$ of 
a finite vector space $({\Bbb F}_p)^n$.
The {\sl rate} $r$ of the code is ${\rm dim}({\cal C})/n$.
The {\sl Hamming distance} $d$ of ${\cal C}$ is the
smallest number of non-zero co-ordinates in a non-trivial
element of ${\cal C}$. One of the main goals of 
coding theory is to construct codes with large rate and
large Hamming distance. Specifically, an infinite collection of codes
is known as {\sl asymptotically good} if $r/n$ and $d/n$ are both bounded
away from zero. The construction of asymptotically good sequences of codes is an
interesting and difficult problem. They were first proved
to exist using probabilistic methods, but explicit
constructions are now available ([2],[10]).

In our situation, the ambient vector space $V$ of the code 
is the space of cellular 1-dimensional mod $p$ cochains on $K_i$. It has
a natural basis, where each basis element is supported on a
single 1-cell. Hence, its dimension is
equal to the number of 1-cells of $K_i$. Pick a
basis for $H^1(K_i; {\Bbb F}_p)$, and represent
each element by a cellular cocycle. The subspace
of $V$ spanned by these cocycles we view as the code
${\cal C}_i$. Let $n_i$ be the dimension of $V$,
and let $r_i$ and $d_i$ be the rate and Hamming
distance of ${\cal C}_i$. The assumption that
$\{ \pi_1(K_i) \}$ has linear growth of mod $p$ homology
is equivalent to the statement that $r_i/n_i$ is
bounded away from zero. The quantity $d_i/n_i$
simply measures the smallest ratio between the support
size of a non-trivial cocycle in ${\cal C}_i$ and
the number of 1-cells of $K_i$. Hence, it is an upper
bound for the smallest relative size of a non-trivial
class in $H^1(K_i ;{\Bbb F}_p)$. Thus, we have the following.

\noindent {\bf Theorem 6.3.} {\sl
Let $K$ be a finite connected
2-complex, and let $\{ K_i \rightarrow K\}$ be a
collection of finite-sheeted covering spaces.
Suppose that $\{ \pi_1(K_i) \}$ has linear
growth of mod $p$ homology for some prime $p$.
Suppose also that there is some $\epsilon > 0$ such that
the relative size of any non-trivial class in $H^1(K_i; {\Bbb F}_p)$
is at least $\epsilon$. Then the codes
${\cal C}_i$ described above are asymptotically good.}

Combining Theorems 6.1 and 6.3, we have the following
result. 

\noindent {\bf Theorem 6.4.} {\sl Let $K$ be a finite
connected 2-complex, and let $\{ K_i \rightarrow K\}$ be a
collection of finite-sheeted covering spaces.
Suppose that $\{ \pi_1(K_i) \}$ has linear
growth of mod $p$ homology for some prime $p$.
Then either $\pi_1(K)$ is large or the codes
${\cal C}_i$ described above are asymptotically good.}

\vfill\eject
\centerline{\caps  7. Finitely generated versus finitely presented}
\vskip 6pt

In Theorem 1.1, we assumed that $G$ was finitely presented.
The remaining hypotheses make sense when $G$ is only
finitely generated. So, it is natural to enquire
whether Theorem 1.1 remains true when the
hypothesis of being finite presented is weakened
to being finitely generated. In this section, we show
that the answer is `no', by analysing a collection
of examples. These were suggested to the author by
Jim Howie. Using the same examples, we also
show that the hypothesis of finite presentability
cannot be weakened in Theorem 1.3 and Corollary 1.8.
The argument here was supplied by Alex Lubotzky.

The groups we will study are the generalised
lamplighter groups $({\Bbb Z}/p{\Bbb Z}) \wr {\Bbb Z}$.
(When $p = 2$, this is the usual lamplighter
group.) Each is a semi-direct product 
$(\oplus_{-\infty}^\infty ({\Bbb Z}/p{\Bbb Z}))
\rtimes {\Bbb Z}$. Here, an arbitrary element
of $\oplus_{-\infty}^\infty ({\Bbb Z}/p{\Bbb Z})$
is required to have only finitely many non-zero co-ordinates.
To define the semi-direct product, we must specify the
action of ${\Bbb Z}$ on $\oplus_{-\infty}^\infty ({\Bbb Z}/p{\Bbb Z})$.
The action of an integer $n$ in ${\Bbb Z}$ 
on $\oplus_{-\infty}^\infty ({\Bbb Z}/p{\Bbb Z})$
simply shifts the indexing set $n$ to the right.
These groups are finitely generated but not finitely
presented [1]. Indeed, each is generated by two elements $a$ and $b$,
where $a$ shifts the indexing set one to the right,
and $b$ lies in $\oplus_{-\infty}^\infty ({\Bbb Z}/p{\Bbb Z})$,
with a single non-zero entry which takes the value
1 in the zero copy of ${\Bbb Z}/p{\Bbb Z}$.

\noindent {\bf Proposition 7.1.} {\sl The generalised lamplighter
group $G = ({\Bbb Z}/p{\Bbb Z}) \wr {\Bbb Z}$ has a nested sequence of finite index
normal subgroups $\{ G_i \}$, each with index a power
of $p$, with the following properties:
\item{(i)} $G$ does not have Property $(\tau)$ with respect to $\{ G_i \}$,
and
\item{(ii)} $\{ G_i \}$ has linear growth of mod $p$ homology.

\noindent But $G$ is not large.}

\noindent {\sl Proof.} By the definition of the
semi-direct product, $G$ admits a surjective
homomorphism $\phi$ onto ${\Bbb Z}$. Let $G_i$
be $\phi^{-1}(p^i {\Bbb Z})$. Then $G_i$
is normal and has index $p^i$. Clearly,
these subgroups are nested.

(i): Lemma 2.1
states that $G$ has property $(\tau)$ with respect
to $G_i$ if and only if ${\Bbb Z}$ has property
$(\tau)$ with respect to $\{ p^i {\Bbb Z} \}$.
But, we have already seen in the example in
Section 2 that this is not the case.

(ii): We claim that $d_p(G_i) \geq [G:G_i]$.
To do this, we will find $p^i$ linearly independent
homomorphisms $G_i \rightarrow {\Bbb F}_p$.
Now, $G_i$ is the subgroup of $G$ generated
by $\oplus_{-\infty}^\infty ({\Bbb Z}/p{\Bbb Z}) $
and $a^{p^i}$. Each homomorphism will send $a^{p^i}$ to
the identity. To define such a 
homomorphism, it suffices to define a homomorphism
$\oplus_{-\infty}^\infty ({\Bbb Z}/p{\Bbb Z}) \rightarrow {\Bbb F}_p$
which is invariant under the action of $a^{p^i}$.
Let $j$ be an integer between $0$ and $p^i-1$.
Define
$$\eqalign{
\oplus_{-\infty}^\infty ({\Bbb Z}/p{\Bbb Z}) &\buildrel
\phi_j \over \longrightarrow {\Bbb F}_p \cr
(n_k)_{k=-\infty}^\infty &\mapsto \sum_{k = -\infty}^\infty n_{{p^i}k + j}.}$$
These are clearly linearly independent, as required. 

Finally, $G$ is not large, because it is soluble. $\square$

We now show that Theorem 1.3 does not remain true
for finitely generated, infinitely presented groups.

\noindent {\bf Proposition 7.2.} {\sl The generalised
lamplighter group $G = ({\Bbb Z}/p{\Bbb Z}) \wr {\Bbb Z}$
does not have Property $(\tau)$ with respect to
any infinite collection of finite index subgroups.}

However, as we have seen in Proposition 7.1, $G$ does
have a nested sequence of normal subgroups, each with
index a power of $p$, that have linear growth of mod
$p$ homology. Hence, by Corollary 1.8, the pro-$p$
completion of $G$ has exponential subgroup growth.

\noindent {\sl Proof of Proposition 7.2.} 
Now, $G$ is amenable, and Theorem 3.1 of [8] asserts
that a finitely generated amenable group does not
have Property $(\tau)$ with respect to any infinite
family of finite index normal subgroups. However,
the assumption of normality is not required in
the proof of that theorem. The proposition now follows.
$\square$

I am grateful to Alex Lubotzky who informed me of his work with Weiss [8],
which formed the basis for this proof.

\vfill\eject
\centerline{\caps 8. Subgroup growth and linear growth of homology}
\vskip 6pt

Throughout this paper, the main focus has been on groups
having a sequence of subnormal subgroups, each with index a power
of a prime $p$, and with linear growth of mod $p$ homology.
In this section, we show how the existence of such a sequence
of subgroups has equivalent characterisations in terms
of subgroup growth.

For a group $G$, let $s_n(G)$ be the number of subgroups with
index at most $n$, and let $a_n(G)$ be the number of
subgroups with index precisely $n$. Let $s_n^{\triangleleft \triangleleft}(G)$
and $a_n^{\triangleleft \triangleleft}(G)$ be the number of
subnormal subgroups with index at most $n$ and precisely
$n$, respectively. A group is said
to have {\sl (at least) exponential subgroup growth} if 
$$\limsup_n {\log s_n(G) \over n} > 0.$$ If $p$ is a prime,
let $\hat G_{(p)}$ be the pro-$p$ completion of $G$.
It turns out that the subgroup growth of a finitely
generated pro-$p$ group is at most exponential.
In other words, $\limsup_n \log s_n(\hat G_{(p)}) / n$ is finite
(Theorem 3.6 of [7]).

The following is a stronger version of Theorem 1.7, which
was stated in the Introduction.

\noindent {\bf Theorem 8.1.} {\sl Let $G$ be a finitely generated
group, and let $p$ be a prime. Then the following
are equivalent:
\item{(i)} $G$ has an infinite nested sequence of subnormal subgroups,
each with index a power of $p$, and with linear growth
of mod $p$ homology;
\item{(ii)} $\hat G_{(p)}$ has exponential subgroup growth;
\item{(iii)} $\limsup_n (\log a_{p^n}^{\triangleleft \triangleleft}(G))/p^n > 0$;
\item{(iv)} $\inf_{n\geq 1} (\log a_{p^n}^{\triangleleft \triangleleft}(G))/p^n > 0$.

}

\noindent {\sl Proof.} (ii)$\Leftrightarrow$(iii): 
There is a one-one correspondence
between subnormal subgroups of $G$ with index $p^n$ and
subgroups of $\hat G_{(p)}$ with index $p^n$. In addition,
any finite index subgroup of $\hat G_{(p)}$ has index
a power of $p$.
Thus, the equivalence of (ii) and (iii) is consequence
of the following general fact. Any sequence of non-negative integers $c_j$
has at least exponential growth (that is,
$\limsup_j (\log c_j)/j > 0$) if and only if the
partial sums $\sum_{i=0}^j c_i$ have at least exponential
growth. In this case, the sequence $c_j$ is
$a_{j}^{\triangleleft \triangleleft}(G)$ if $j$ is
a power of $p$, and zero otherwise.

(i)$\Rightarrow$(iv):
Suppose that $G$ has a sequence of subgroups
$G = G_1 \triangleright G_2 \triangleright \dots$
such that $G_n/G_{n+1}$ is a non-trivial finite $p$-group for each $n$,
and with linear growth of mod $p$ homology.
Let $\lambda$ be $\inf_n (d_p(G_n) -1)/[G:G_n]$, 
the mod $p$ homology gradient, which 
is therefore positive.
Now, any finite $p$-group has a subnormal series,
where successive quotients are cyclic of order $p$.
Thus, by refining the sequence $\{ G_n \}$ if necessary, we may
assume that each $G_n/G_{n+1}$ is cyclic of order
$p$. By Proposition 5.1, $(d_p(G_n) - 1)/[G:G_n]$ is
a non-increasing function of $n$. Thus, 
$\inf_n (d_p(G_n) -1)/[G:G_n]$ is still $\lambda$. Any normal
subgroup of $G_n$ with index $p$ arises as the kernel
of a non-trivial homomorphism $G_n \rightarrow {\Bbb Z}/p{\Bbb Z}$.
There are $p^{d_p(G_n)} -1$ such homomorphisms. The number
of homomorphisms with a given kernel is $p-1$. Thus,
there are $(p^{d_p(G_n)} -1)/(p-1)$ normal subgroups of
$G_n$ with index $p$. Each gives
a subnormal subgroup of $G$ with index $[G:G_n]p = p^n$.
Hence, when $n \geq 1$,
$a_{p^n}^{\triangleleft \triangleleft}(G)$ is
at least 
$${p^{\lambda p^{n-1}+1} - 1 \over p-1},$$ and so we deduce 
that $\liminf_n (\log a_{p^n}^{\triangleleft \triangleleft}(G))/p^n$
is positive. Finally, note that
$a_{p^n}^{\triangleleft \triangleleft}(G)$
is always more than 1, when $n \geq 1$, and so
$(\log a_{p^n}^{\triangleleft \triangleleft}(G))/p^n$
is strictly positive. Thus, we deduce (iv).

(iv)$\Rightarrow$(iii): This is trivial. 

(iii)$\Rightarrow$(i): Define
$$r_n = \max \{ d_p(H): H \triangleleft\!\triangleleft \ G
\hbox{ and } [G:H] = p^n \}.$$
Let us suppose that (i) does not hold.
We claim that $\limsup_n r_n / p^n = 0$. 
For otherwise, $\limsup_n r_n / p^n$ is positive, and
therefore so is $\limsup_n (r_n-1) / p^n$. Let
$\lambda$ be this latter value. Note that,
by Proposition 5.1, $(r_n - 1)/p^n$ is a
non-increasing function of $n$. Thus,
$\lambda$ is actually the infimum and limit
of this sequence.
Hence, for each $n$, there is a subnormal
subgroup $G_n$, with index $p^n$ such that
$d_p(G_n) -1\geq \lambda p^n$. For each $n$,
we may find a subnormal sequence
$$G = G_{n,1} \triangleright G_{n,2} \triangleright
\dots \triangleright G_{n,n} = G_n$$ such that
$G_{n,i}/G_{n,i+1}$ is cyclic of order $p$ for each $i$.
Now, $(d_p(G_{n,i}) - 1) /p^i \geq \lambda$ by
Proposition 5.1. Since $G$ has only finitely many
subgroups of index $p$, we may find a subsequence of
the $G_n$ where $G_{n,2}$ is a fixed group $G_2$. By passing
to a further subsequence, we may assume that $G_{n,3}$
is a fixed group $G_3$, and so on. Thus, we obtain a
sequence of subnormal subgroups $G = G_1 \triangleright G_2 \triangleright \dots$, each with
index $p$ in its predecessor, and with linear growth of
mod $p$ homology. This is condition (i), which we are
assuming does not hold. This contradiction proves the
claim: $\limsup_n r_n / p^n = 0$. Hence, $\lim_{n \rightarrow \infty} \left( \sum_{i=0}^n r_i \right) / p^n = 0$.
Now, any subnormal subgroup of $G$ with index $p^n$
is a normal subgroup of some subnormal subgroup of $G$
with index $p^{n-1}$. Hence, 
$$a_{p^n}^{\triangleleft \triangleleft}(G)
\leq a_{p^{n-1}}^{\triangleleft \triangleleft}(G) p^{r_{n-1}}.$$
Thus, by induction,
$$a_{p^n}^{\triangleleft \triangleleft}(G)
\leq p^{\sum_{i=0}^{n-1} r_i}.$$
Taking logs:
$$\log a_{p^n}^{\triangleleft \triangleleft}(G)
\leq (\log p ){\sum_{i=0}^{n-1} r_i}.$$
Therefore, $$\left( \log a_{p^n}^{\triangleleft \triangleleft}(G) \right)  /p^n \rightarrow 0,$$
which means that (iii) does not hold, as required.
$\square$

\vskip 18pt
\centerline{\caps References}
\vskip 6pt

\item{1.} {\caps G. Baumslag,}
{\sl Wreath products and finitely presented groups},
Math. Z. 75 (1960/1961) 22--28.

\item{2.} {\caps J. Justesen}, {\sl A class of constructive
asymptotically good algebraic codes}, IEEE Trans. Inf. Theory
IT-18 (1972) 652--656.

\item{3.} {\caps M. Lackenby}, {\sl Heegaard splittings,
the virtually Haken conjecture and Property $(\tau)$},
Invent. Math. (to appear)

\item{4.} {\caps M. Lackenby}, {\sl Expanders, rank and
graphs of groups}, Israel J. Math. 146 (2005) 357--370.

\item{5.} {\caps M. Lackenby}, {\sl Covering spaces of 3-orbifolds},
Preprint.

\item{6.} {\caps A. Lubotzky}, {\sl Discrete Groups,
Expanding Graphs and Invariant Measures}, Progress
in Math. 125 (1994)

\item{7.} {\caps A. Lubotzky and D. Segal}, {\sl Subgroup growth},
Progr. in Math. 212 (2003)

\item{8.} {\caps A. Lubotzky and B. Weiss}, {\sl Groups and
expanders},  Expanding graphs (Princeton, NJ, 1992),  95--109,
DIMACS Ser. Discrete Math. Theoret. Comput. Sci., 10.

\item{9.} {\caps A. Lubotzky and R. Zimmer}, {\sl
Variants of Kazhdan's property for subgroups of semisimple groups}, 
Israel J. Math.  66 (1989) 289--299.

\item{10.} {\caps D. Spielman,} {\sl Constructing error-correcting codes
from expander graphs},  Emerging applications of number theory 
(Minneapolis, MN, 1996),  591--600,
IMA Vol. Math. Appl., 109.

\vskip 12pt
\+ Mathematical Institute, University of Oxford, \cr
\+ 24-29 St Giles', Oxford OX1 3LB, United Kingdom. \cr

\end